\documentclass{amsart}
\usepackage{amsfonts,amsmath,amssymb,amscd,latexsym,graphicx,epsfig,color}

\newtheorem{thm}{Theorem}[section]

\newtheorem{cor}[thm]{Corollary}

\newtheorem{example}[thm]{Example}

\newtheorem{rmk}[thm]{Remark}

\def\Z{\mathbb{Z}}
\def\R{\mathbb{R}}
\def\N{\mathbb{N}}

\def\co{\colon\thinspace}

\newcommand{\esp}{\hspace{0.1cm}}

\begin{document}

\title[Nielsen-Thurston orderings and the space of braid orderings]{Nielsen-Thurston orderings and\\ the space of braid orderings}

\author{Andr\'es Navas}
\address{Andr\'es Navas, 
Univ. de Santiago de Chile, Alameda 3363, Santiago, Chile}
\email{andres.navas@usach.cl}

\author{Bert Wiest}
\address{Bert Wiest, IRMAR (UMR 6625 du CNRS), Universit\'e de Rennes 1,
Campus de Beaulieu, 35042 Rennes Cedex, France}
\email{bertold.wiest@univ-rennes1.fr}




\begin{abstract}
We study the topological space of left-orderings of the braid group,
and its subspace of Nielsen-Thurston orderings. Our main result is
that no Nielsen-Thurston ordering is isolated in the space of braid
orderings. In the course of the proof, we classify the convex subgroups
and calculate the Conradian soul for any Nielsen-Thurston ordering 
of~$B_n$.
We also prove that for a large class of Nielsen-Thurston orderings, 
including all those of infinite type, a stronger result holds: they
are approximated by their own conjugates. 
On the other hand, we suggest an example of a Nielsen-Thurston ordering
which may not be approximated by its conjugates.
\end{abstract}
\maketitle


\section{Introduction}

In recent years, some progress in the theory of orderable groups has been achieved 
by studying the so-called {\em space of orderings} (throughout, the word 
{\em ordering} stands for a left-invariant total order relation on a group). 
This corresponds to the set of all orderings that may be defined on a 
given group $\Gamma$, and carries a natural topology: given finitely 
many group elements, a neighborhood of a prescribed ordering $\prec$ 
is the set of all orderings which coincide with $\prec$ on this finite set. 
It is easy to see that this space is totally disconnected and compact; 
moreover, it is metrizable when the underlying group is countable. 
It is a  nontrivial fact that this space cannot be countably infinite 
\cite{linnell,navas-dynamics,CNR}.
Notice that an isolated point in this space corresponds to an ordering 
which is completely determined by finitely many inequalities.

Remarkable examples of orderable groups are the braid groups $B_n$. 
Historically, the first ordering on $B_n$ (for $n \geqslant 3$) was defined 
by Dehornoy. It was  shown in~\cite{navas-dynamics} that the Dehornoy 
ordering $\prec_{_D}$ is non isolated in the space of braid orderings. 
This is to be contrasted with the work~\cite{DD}, where Dubrovina and 
Dubrovin show, by means of a slight modification of the Dehornoy 
construction, that $B_n$ admits isolated orderings. For example, 
on $B_3$ there is a unique ordering $\prec_{_{D\hspace{-1pt}D}}$ 
verifying the inequalities \esp 
$\sigma_1 \sigma_2 \succ\hspace{-3pt}_{_{D\hspace{-1pt}D}} 1$ \esp and 
\esp $\sigma_2^{-1} \succ\hspace{-3pt}_{_{D\hspace{-1pt}D}} 1.$

Dehornoy's construction involves some deep ideas which are purely 
algebraic. However, 
there is an alternative geometric approach to braid orderability. 
Indeed, $\prec_{_D}$ can be seen as a particular member of a family of 
uncountably many orderings, each of which is associated to a completely 
separating geodesic of the punctured disk (endowed with a hyperbolic 
metric). These so-called {\em Nielsen-Thurston 
orderings} on braid groups where largely studied and classified by 
H. Short and the second-named author in~\cite{SW}. In the present 
paper we pursue this study from the point of view of the space of 
orderings. Our main result may be stated as follows.

\medskip

\noindent{\bf Theorem.} {\em No Nielsen-Thurston ordering is isolated in the space 
of braid-orderings.}

\medskip

We recall that the \emph{positive cone} of an ordering is the
semigroup consisting of elements greater than the identity.
As a direct consequence of this theorem and \cite[Proposition 1.8]{navas-dynamics} we have:

\medskip

\noindent{\bf Corollary.} {\em 
There is no Nielsen-Thurston ordering whose positive cone is finitely 
generated as a semigroup.}

\medskip

To show the theorem, we first determine the convex subgroups of Nielsen-Thurston 
orderings (where a subset $S$ is said to be {\em convex} if $g \in S$ whenever 
$f_1 \prec g \prec f_2$ and both $f_1,f_2$ are in $S$). We prove that the only convex 
subgroups are the obvious ones, coming from the way in which the geodesic cuts up 
the surface. For example, for the Nielsen-Thurston orderings of 
{\em full} infinite type (to be defined below), we prove that they 
have no convex subgroups other than $\{1\}$ and $B_n$. This 
solves by the affirmative the second half of Question 2.21 in \cite[Chapter XVI]{DDRW}. 

With the whole list of convex subgroups at hand, for each Nielsen-Thurston 
ordering we determine the {\em Conradian soul}, that is the maximal convex 
subgroup restricted to which the ordering satisfies the so-called 
{\em Conrad property}. This notion was introduced by the first-named author 
as a tool for studying the possibility of approximating a given ordering 
by its conjugates (where the conjugate of an ordering $\prec$ by a group 
element $h$ is the ordering $\prec_h$ for which $f \prec_h g$ if and 
only if $h^{-1} f h \prec h^{-1} g h$). 
For example, by~\cite{navas-dynamics} every ordering on a countable
group having trivial 
Conradian soul may be approximated by its conjugates. 
As a consequence, Nielsen-Thurston orderings of full 
infinite type can be approximated by their conjugates.

For Nielsen-Thurston orderings whose Conradian soul is isomorphic 
to~$\mathbb{Z}$ (as is the case of the Dehornoy ordering), we show 
that the property of accumulation by its conjugates still holds. 
In the rest of the cases, we show that the Conradian soul is  
isomorphic to $\mathbb{Z}^k$ for some $k \geqslant 2$. Although 
we ignore whether these orderings may be approximated by their 
conjugates (actually, we present evidence that some of them might be
isolated in the ``space of Nielsen-Thurston orderings''), we show 
that they are not isolated in the space of braid orderings. This 
is achieved  by means of a classical and very simple 
{\em convex extension} type argument.


\section{A brief reminder on Nielsen-Thurston orderings}

Nielsen-Thurston orderings of $B_n$, which were introduced in~\cite{SW},
are the total orderings arising from the natural action of $B_n$ on the
open interval $]0,\pi[$. 

We briefly recall the definition. Let $D_n$ denote the unit disk, 
with~$n$ punctures lined up on the horizontal diameter.
We can equip~$D_n$ with a hyperbolic 
structure (in which the punctures correspond to cusps). Indeed, many 
such hyperbolic structures exist, but we fix one choice, for the 
rest of the paper. 

A Nielsen-Thurston ordering is, by definition, induced by a 
finite or infinite geodesic 
ray~$\gamma\co [0,T]\to D_n$ (with $T\in\R_+$) or 
$\gamma\co [0,\infty[\to D_n$, starting at a basepoint on the boundary 
of~$D_n$, and (in the case of a finite ray) terminating at some other 
point of~$\partial D_n$.

The ordering induced by~$\gamma$ is defined as follows. 
Any geodesic starting on the basepoint in $\partial D_n$ has in its 
starting point a certain angle, belonging to the interval $]0,\pi[$,
with the boundary $\partial D_n$. Conversely, the geodesic is uniquely
determined by the value of this angle. More generally, even a
non-geodesic ray starting on the basepoint is homotopic, relative
to its end points, to a unique geodesic, so homotopy classes of
such rays are still in one-to-one correspondence with angles in
$]0,\pi[$. Now the braid group, viewed as the
mapping class group of~$D_n$, acts on the geodesic~$\gamma$; 
this yields a partial ordering $\prec$ of $B_n$, by defining that 
$\beta_1\prec \beta_2$ if and only if the angle of $\beta_1.\gamma$
is smaller than the angle $\beta_2.\gamma$. 

If $\gamma$ is not invariant under the action of any braid, then 
this recipe determines a \emph{total} ordering of~$B_n$, and this 
ordering is called \emph{the Nielsen-Thurston ordering induced by~$\gamma$}.
Notice that a necessary (but not sufficient) condition for~$\gamma$
not being invariant under any braid action is 
that~$\gamma$ cuts up the surface, in the sense that no two punctures are 
in the same path components of~$D_n-\gamma([0,T])$). 

An important property of Nielsen-Thurston orderings is the so-called 
{\em subword-property}, or {\em property~S}. This means that for 
any two braid words $w$ and $w'$, where~$w'$ is obtained from~$w$ by 
inserting positive generators~$\sigma_i$, we have~: $w \prec w'$. 
Geometrically, this is nothing but the fact that the left half-twist 
along any arc connecting two punctures pushes any curve intersecting
the arc further to the left. In particular this half-twist is a positive 
element for the ordering.   

The properties of the Nielsen-Thurston orderings are closely related to 
the geometric nature of the corresponding geodesic. Let us consider, 
for a maximal $m\in\{0,1,\ldots,n-1\}$, the~$m$ points in time 
such that some pair of 
punctures is in the same path component of~$D_n-\gamma([0,t_i-\epsilon])$, 
but in separate path components of~$D_n-\gamma([0,t_i])$. We shall call 
this the \emph{sequence of separating moments}.

If $m=n-1$, i.e.~if a finite initial segment of~$\gamma$ is 
sufficient to cut up the surface, then the ordering is said to be of 
\emph{finite type}. If $m<n-1$ then the ordering is of \emph{infinite
type}. If $m=0$, i.e.~if $\gamma$ has no self-intersections and no finite
initial segment of~$\gamma$ separates any pair of punctures, then we say
the ordering is of \emph{full infinite type}.


\section{Convex subgroups of Nielsen-Thurston orderings of finite type}

\begin{thm}\label{T:ConvFinite}
Suppose $\prec$ is a Nielsen-Thurston ordering of finite type, associated to a
geodesic~$\gamma$. Let $0<t_1<t_2<\ldots<t_{n-1}$ be the sequence of
separating moments of~$\gamma$.
Then the chain of convex subgroups is
$$
\{1\}\subset G_{n-1}\subset\ldots\subset G_1\subset B_n
$$
where $G_i$ denotes the subgroup of~$B_n$ consisting of all braids
which preserve the geodesic segment $\gamma([0,t_i])$. 
\end{thm}

{\bf Proof.} It is immediate that the subgroups thus described are convex,
we only have to prove that they are the only ones. 

The proof is by downwards induction.  The start of the induction is
with the cases where all components of~$D_n-\gamma([0,t_{n-i}])$ contain
at most two punctures. In these cases, the subgroup consisting of
elements of~$B_n$ leaving~$\gamma_{n-i}$ invariant is isomorphic to~$\Z^i$,
and is lexicographically ordered. Such an ordered group has no convex
subgroups except the ones obtained by successively removing the
largest generator -- and these subgroups do indeed appear in the
chain of convex subgroups $G_{n-1}\subset\ldots\subset G_{n-i}$.

Let us now suppose inductively that for some index~$i$ there are
no convex subgroups in~$G_{n-i}$ except~$G_{n-j}$ with~$j<i$. 
Let $G'$ be a convex subgroup of $G_{n-i-1}$ such that 
$G_{n-i}\subsetneq G'\subset G_{n-i-1}$. 
Our aim is to prove that $G'=G_{n-i-1}$.

Let us recall the structure of the group~$G_{n-i-1}$: the 
geodesic~$\gamma([0,t_{n-i-1}])$ cuts~$D_n$ into a number of components,
some of which contain punctures. If we denote $n_1,\ldots,n_k$ the
number of punctures contained in each of these path components, then the 
subgroup~$G_{n-i-1}$ is isomorphic to a product
$G_{n-i-1}\cong B_{n_1}\times\ldots\times B_{n_k}$. Moreover, there
is exactly one connected component of $D_n - \gamma([0,t_{n-i-1}])$
which is further cut by $\gamma([t_{n-i-1},t_{n-i}])$ into pieces
one of which contains~$\lambda$ punctures and another one~$\mu$ 
punctures (with $\lambda,\mu\geqslant 1$, and at least one of the
two being $\geqslant 2$). Without loss of generality let us
say the component being cut is the first one, so that $n_1=\lambda+\mu$.

\begin{figure}[htb]
\centerline{\input{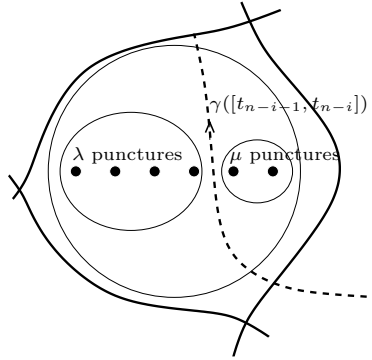}}
\caption{Fat lines: $\gamma([0,t_{n-i-1}])$. Dashed line: 
$\gamma([t_{n-i-1},t_{n-i}])$.}
\label{F:CriticalRegion}
\end{figure}

Now there are three disks with geodesic boundary contained in
the region under discussion: one containing~$\lambda$ punctures, one 
containing~$\mu$ punctures, and one containing both sets of punctures.
Possibly after a conjugation we can
assume that the three disks are \emph{round} 
(i.e. their boundaries intersect the horizontal diameter only twice, see 
Figure~\ref{F:CriticalRegion}).

We shall study the action of the hypothetical subgroup~$G'$ on the 
geodesic~$\gamma([0,t_{n-i}])$, keeping in mind that it fixes the
initial segment~$\gamma([0,t_{n-i-1}])$. 
Let 
$$\gamma_{max}=\sup_{g'\in G'} g'.\gamma$$ 
(where the supremum is in the sense that every geodesic corresponds
to a point of the real line, and we're taking the supremum of real
numbers). What does the first intersection of $\gamma_{max}$ with our
critical component of $D_n- \gamma([0,t_{n-i-1}])$ look like?
Notice that $\gamma_{max}$ must be $G'$-invariant.

\begin{figure}[htb]
\centerline{\input{GeodInCritRegion.pstex_t}}
\caption{}
\label{F:GeodInCritRegion}
\end{figure}

The geodesic segment $\gamma_{max}$ cannot pass between two punctures
among the $\lambda$ punctures of the left disk, nor can it pass between
two of the $\mu$ punctures of the right disk, for if it did, then the
action of an appropriate element in the smaller convex subgroup~$B_\lambda$
or $B_\mu$ would displace the geodesic even further to the left, which
is excluded by construction (see Figure~\ref{F:GeodInCritRegion}(a)). 
Suppose, on the other hand, that $\gamma_{max}$ avoids entirely the 
two discs with geodesic boundary containing the set of~$\lambda$ or $\mu$ 
punctures, but intersects the horizontal segment between the two
disks. Then we have 
$g'.\gamma_{max}>\gamma_{max}$ or $(g')^{-1}.\gamma_{max}>\gamma_{max}$, 
where $g'$ is any element of $G'$ which moves the 
segment~$\gamma([t_{n-i-1},t_{n-i}])$ (such an element must 
exist, since $G_{n-i} \subsetneq G'$). This is also impossible.

In summary, the geodesic $\gamma_{max}$ remains disjoint from the disk 
with geodesic boundary containing the $\lambda+\mu$ punctures; thus it 
spirals onto its boundary, as indicated in Figure~\ref{F:GeodInCritRegion}(b). 
But this means that $g'.\gamma$ can get as big (far to the left) as any 
geodesic $g.\gamma$ with $g\in G_{n-i-1}$. This implies that $G'= G_{n-i-1}$.
\hfill $\Box$


\section{Convex subgroups of Nielsen-Thurston orderings of full infinite type}

We recall that a geodesic ray~$\gamma$ in~$D_n$ gives rise to an ordering
of full infinite type if it starts at the basepoint in~$\partial D_n$, 
has no self-intersections, separates the punctures, and is not stabilized 
by any nontrivial element of~$B_n$.

The next result should be compared with \cite{clay2}.

\medskip

\begin{thm}\label{T:ConvInfinite}
Suppose~$\gamma$ is an infinite geodesic ray in~$D_n$ giving rise to a
Nielsen-Thurston ordering~$\prec$ of~$B_n$ of full infinite type. Then 
$\prec$ has no convex subgroups except $\{1\}$ and $B_n$.
\end{thm}

{\bf Proof.} The strategy is to prove that~$B_n$ can be generated by
elements which are ``arbitrarily close to the identity'', which implies
that any nontrivial convex subgroup contains all of~$B_n$. 

Here are the details. The geodesic~$\gamma$ can be specified by an infinite
word, where the $k$th letter ($k\in\N$) specifies between which
pair of punctures of~$D_n$ the $k$th intersection of~$\gamma$ with the
horizontal diameter occurs, and in which direction (up or down). 

We shall prove that the braid group can be generated by~$n-1$ elements
$\tau_1,\ldots,\tau_{n-1}$ which all leave an arbitrarily long initial 
segment of~$\gamma$ invariant: the words describing $\tau_i.\gamma$
all coincide with the word describing $\gamma$ on an arbitrarily long
initial segment.

Here is the construction of the elements~$\tau_i$: we trace out the 
geodesic~$\gamma$, for at least the whole required initial segment, and
continue until we pass close to a puncture (so close that the line
we have drawn so far does not intersect the short segment from our
current position to the puncture). At that moment, we deviate from the
geodesic~$\gamma$, and drop into the puncture. We replace the curve we
have just drawn by a geodesic which is homotopic to it relative to its 
endpoints, and call this geodesic~$\gamma'$. The puncture at its
end will be called the \emph{central puncture}. 
Let~$k$ be the
number of intersections of~$\gamma'$ with the horizontal diameter
(not counting the final drop into the central puncture).

There are geodesic arcs $a_1,\ldots,a_{n-1}$ disjoint from $\gamma'$
and from each other, connecting the central puncture to the $n-1$
other punctures. For later use, we shall equip these arcs with an
orientation pointing towards the central puncture.
The braid~$\tau_i$ will be a positive or negative
half Dehn-twist around the arc~$a_i$ -- the sign remains to be
specified, depending on the relative position of terminal segments
of~$\gamma'$ and~$a_i$. 

We shall consider three cases.
Firstly, if the last intersections of these two arcs with the horizontal
diameter do \emph{not} occur between the same pairs of punctures,
then the sign of the half Dehn twist can be chosen arbitrarily.
Indeed, in this case the first~$k$ intersections of $\tau_i.\gamma$
with the horizontal diameter coincide with those of~$\gamma$,
independently of whether~$\tau_i$ is taken to be the positive or
negative half twist along~$a_i$ -- see Figure~\ref{F:TwistAction}(a).

\begin{figure}[htb]
\centerline{\input{TwistAction.pstex_t}}
\caption{}
\label{F:TwistAction}
\end{figure}

The second case to consider is where the last intersections of the 
the arcs~$\gamma'$ and~$a_i$ with the horizontal diameter lie between
the same pairs of punctures, and where, moreover, the terminal segment
of~$a_i$ lies to the right of the terminal segment of~$\gamma'$ (as seen
while approaching the central puncture along~$\gamma'$)  -- see 
Figure~\ref{F:TwistAction}(b) for examples. In this case, we take~$\tau_i$ 
to be the positive half twist along~$a_i$. Again, we observe that with
this choice the word describing the first~$k$ intersections of 
$\tau_i.\gamma$ with the horizontal diameter coincides with the
corresponding word for~$\gamma$. 

The third and final case is where the last intersections of the 
the arcs~$\gamma'$ and~$a_i$ with the horizontal diameter lie between
the same pairs of punctures, and where the terminal segment
of~$a_i$ lies to the \emph{left} of the terminal segment of~$\gamma'$.
In this case, we take~$\tau_i$ to be the negative half twist along~$a_i$
-- see Figure~~\ref{F:TwistAction}(c).
\hfill$\Box$
\smallskip


\section{Convex subgroups of general Nielsen-Thurston orderings}

By combining the techniques of proof of Theorems~\ref{T:ConvFinite} 
and~\ref{T:ConvInfinite}, one can easily obtain the following 
unified result:

\medskip

\begin{thm}\label{T:ConvGeneral}
Suppose $\prec$ is a Nielsen-Thurston ordering, associated to a
geodesic~$\gamma$. Let $0<t_1<t_2<\ldots<t_m$, with
$m\in\{0,1,\ldots,n-1\}$, be the sequence of
separating moments of~$\gamma$. 
Then the chain of convex subgroups is
$$
\{1\}\subset G_m\subset\ldots\subset G_1\subset B_n
$$
where $G_i$ denotes the subgroup of~$B_n$ consisting of all braids
which preserve the geodesic segment $\gamma([0,t_i])$. 
\hfill$\Box$
\end{thm}


\section{The Conradian soul of Nielsen-Thurston orderings}

Recall that a group ordering $\prec$ is said to be Conradian if for every 
positive $f,g$ there exists $k \in \mathbb{N}$ so that $f g^k  \succ g$ 
(see \cite{rhem,kopi}). 
Using structure theorems for Conrad orderable groups and the description 
of the convex subgroups from the previous sections, one can easily show 
that no Nielsen-Thurston ordering is Conradian 
(c.f. Corollary~\ref{no-Conrad}). 
However, we prefer giving a longer proof which avoids these  
structure theorems. In this proof, the cases of finite type orderings on $B_3$ 
and $B_4$ need a special argument (compare Remark \ref{loc-ind}). 

\begin{example}\rm The Dehornoy ordering $\prec_{_D}$ on $B_3$ is not 
Conradian. This is shown in \cite[Example 3.21]{navas-dynamics}, through 
a very indirect argument. A shorter proof works as follows: using the 
identity \esp $\sigma_2 \sigma_1 \sigma_2^{-1} = 
\sigma_1^{-1} \sigma_2 \sigma_1$ \esp  we get, for all $k\in \N$, that
$1 \prec_{_D} \sigma_2 \sigma_1 \sigma_2^{-(k+1)} = 
\sigma_1^{-(k+1)} \sigma_2 \sigma_1$.
Therefore, $(\sigma_2^{-1} \sigma_1) (\sigma_1^k) \prec_{_D} \sigma_1$ for 
all $k \in \mathbb{N}$, yet both $\sigma_2^{-1} \sigma_1$ and $\sigma_1$ 
are $\prec_{_D}$-positive.

Notice that, according to \cite{SW}, every Nielsen-Thurston ordering of 
finite type on $B_3$ 
is conjugate to the Dehornoy ordering.  Since the Conrad property is stable under conjugacy, 
no Nielsen-Thurston ordering of finite type on $B_3$ is Conradian.
\label{one}
\end{example}

\begin{example}\rm Up to conjugacy, there are only three Nielsen-Thurston 
orderings of finite type on $B_4$, namely those given by the geodesics of 
Figure~\ref{F:D4geod} below (see again \cite{SW}). 
None of them gives rise to a Conradian ordering. 
Indeed, (a) corresponds to the Dehornoy ordering, and the above 
argument applies. For the geodesics in (b) and (c), one can directly check 
the inequality \esp $1 \prec \sigma_3 \sigma_2 \sigma_3^{-(k+1)}$, 
\esp which allows to conclude as in the previous example.

\begin{figure}[htb]
\centerline{\input{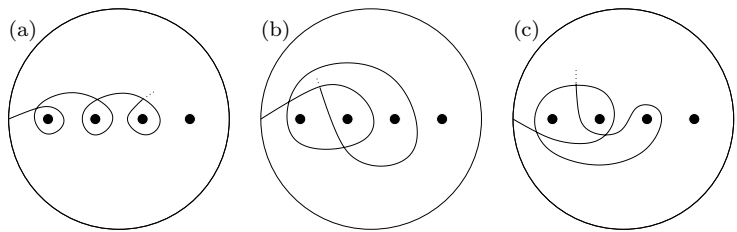}}
\caption{The three conjugacy classes of finite type Nielsen-Thurston orderings of $B_4$.}
\label{F:D4geod}
\end{figure}
\label{two}
\end{example}

Following \cite{navas-dynamics,CNR}, the {\em Conradian soul} of an ordering is defined as 
the maximal convex subgroup restricted to which the order satisfies the Conrad property. 
It was already shown in \cite{navas-dynamics} that the Conradian soul of the 
Dehornoy ordering on $B_n$ corresponds to $\langle \sigma_{n-1} \rangle$. 
Next, we describe the Conradian soul of any Nielsen-Thurston ordering, 
thus solving a problem raised in \cite[Remark 3.40]{navas-dynamics}.

\smallskip

\begin{thm}
The Conradian soul of a Nielsen-Thurston ordering on $B_n$ corresponds to 
the maximal Abelian convex subgroup, that is, to the largest copy of 
$B_2 \times \cdots \times B_2$ which is convex.
\label{descr-CS}
\end{thm}

{\bf Proof.} For a geodesic of infinite type, the smallest non-trivial 
convex subgroup is a copy of $B_m$ for some $m$ with 
$3\leqslant m\leqslant n$. Such a group has no convex subgroup. If it 
was Conradian, then it would be Archimedean, and hence by H\"older's 
theorem, the underlying group would be Abelian (see \cite{rhem,kopi}). 
This is absurd as $B_m$, for $m\geqslant 3$, is not Abelian.

If the geodesic is of finite type, then it is clear that the maximal convex 
copy of a product of (say, $k$) $B_2$ factors is contained in the Conradian 
soul (indeed, any bi-invariant ordering -- and thus any ordering defined on 
an Abelian group -- is Conradian). 
The next larger convex subgroup of $B_n$ corresponds either to a product 
$P = B_{n_1} \times \cdots \times B_{n_{k-1}}$, where there is only one 
index different from $2$, which is equal to $4$, or to a product 
$P = B_{n_1} \times \cdots \times B_{n_k}$, where only one index is different 
from $2$, and it is equal to $3$. If the restriction of the ordering to $P$ 
was Conradian, then the restriction to $B_4$ or $B_3$, respectively, would 
be also Conradian. However, we have seen in Examples \ref{one} 
and \ref{two} that this cannot be the case.

{\em An alternative argument for the end of the proof.} To see that the 
restriction of the ordering to $P$ cannot be Conradian, first note 
that~$P$ has only finitely many convex subgroups. 
By \cite[Proposition 1.7]{CNR} and its proof, if $P$ was 
Conradian then the underlying group would be solvable. 
However, neither $B_3$ nor $B_4$ are solvable. \phantom{i}
\hfill$\Box$
\medskip

\begin{cor} No Nielsen-Thurston ordering is Conradian.
\label{no-Conrad}
\end{cor}

{\bf Proof.} If a Nielsen-Thurston ordering was Conradian, then it would 
coincide with its Conradian soul. However, according to Theorem 
\ref{descr-CS}, this is never the case.
\hfill$\Box$

\smallskip

\begin{rmk}\rm Although stated for Nielsen-Thurston orderings on general braid groups $B_n$ (with $n \geqslant 3$), the preceding corollary 
gives new information only in the cases of $B_3$ and $B_4$. Indeed, 
these groups do admit Conradian orderings, whereas for $n \geqslant 5$ 
there is no Conradian ordering on $B_n$ at all 
(see \cite[pp.~287-289]{DDRW}).  
\label{loc-ind}
\end{rmk}


\section{Nielsen-Thurston orderings of infinite type are 
approximated by their conjugates}

The relevance of the Conradian soul stems from the following 
result which we cite from~\cite[Proposition 4.7]{navas-dynamics} 
(see also \cite[Theorem 1.2]{clay} and \cite[Theorem 2.7]{CNR}).

\medskip

\begin{thm}
If the Conradian soul of an ordering on an infinite countable group 
is trivial, then this ordering is an accumulation point of its conjugates.
\label{CS-thm}
\end{thm}

\medskip

Using this result we may easily prove the following

\medskip

\begin{thm}
Any Nielsen-Thurston ordering of infinite type 
on~$B_n$ may be approximated by its conjugates.
\label{approx-conj}
\end{thm}

{\bf Proof.} By Theorem \ref{descr-CS}, the Conradian soul of the 
ordering is trivial. The result then follows from Theorem \ref{CS-thm}.
\hfill$\Box$


\section{Nielsen-Thurston orderings of finite type are not isolated }

We now show our main result. 

\medskip

\begin{thm}\label{T:NTnotIsol}
Nielsen-Thurston orderings are not isolated in the 
topological space of left-invariant orderings of~$B_n$.
\end{thm}

{\bf Proof.} In the case of orderings of infinite type, the ordering 
may be approximated by its conjugates (see Theorem \ref{approx-conj}).

In the finite type case, for a maximal value of $k$ there exists a convex 
subgroup of the form $B_2 \times \cdots \times B_2 \sim \mathbb{Z}^k$. 
If $k \geqslant 2$, then the restriction of $\prec$ to $ \mathbb{Z}^k$ 
is not isolated in the space of orderings of $\mathbb{Z}^k$ 
(see for instance \cite{sikora}). A convex extension argument (see 
\cite[Section 3.3.5]{navas-dynamics}) then shows that $\prec$ 
is not isolated in the space of braid orderings. If $k = 1$, then the 
smallest convex subgroup strictly containing the maximal Abelian convex 
subgroup $B_2$ is a copy of $B_3$. By \cite{SW}, the 
restriction of $\prec$ to this copy of $B_3$ is a conjugate of the 
Dehornoy ordering, which is not isolated in the space of orderings 
of $B_3$: see Example \ref{cle} below. Once again, a 
convex extension argument shows that $\prec$ is not isolated in 
the space of braid orderings. 
\hfill$\Box$

\medskip

\begin{example}\rm Dehornoy's ordering is approximated by its conjugates.  
This was first shown in \cite{navas-dynamics} throughout a very indirect 
argument. A simpler proof appears in \cite[Chapter XIV]{DDRW}. Here 
we propose an even simpler argument. 

Let $\prec_j$ be the Dehornoy ordering, conjugated by 
$\sigma_2^{-j} \sigma_1$. Thus, a word $w$ is positive in 
the ordering $\prec_j$ if and only if 
$\sigma_1^{-1}\sigma_2^{j}w\sigma_2^{-j}\sigma_1\succ_{_D} 1$.
We claim that the sequence $\prec_j$ tends to $\prec_{_D}$ in the space of orderings.

Indeed, if $w=\sigma_2^k$ for some $k>0$, then 
$\sigma_1^{-1}\sigma_2^{j}w\sigma_2^{-j}\sigma_1 \succ_{_D} 1$ since the Dehornoy
order has Property S. If, on the other hand, $w$ is a $\sigma_1$-positive word:
$w=\sigma_2^{k_1}\sigma_1\sigma_2^{k_2}\sigma_1\ldots\sigma_2^{k_{\ell-1}}\sigma_1
\sigma_2^{k_\ell}$, then we calculate
\begin{eqnarray*}
\sigma_1^{-1}\sigma_2^j w \sigma_2^{-j}\sigma_1 & = &
\sigma_1^{-1}\sigma_2^j \sigma_2^{k_1}\sigma_1\sigma_2^{k_2}\sigma_1\ldots
\sigma_2^{k_{\ell-1}}\sigma_1\sigma_2^{k_\ell} \sigma_2^{-i}\sigma_1\\
 & = &
\sigma_2\sigma_1^{j+k_1}\sigma_2^{-1}\sigma_2^{k_2}\sigma_1\ldots
\sigma_2^{k_{\ell-1}}\sigma_1\sigma_2^{k_\ell} \sigma_2^{-n}\sigma_1.
\end{eqnarray*}
Thus $\sigma_1\sigma_2^{-j}w\sigma_2^j\sigma_1$ is
$\sigma_1$-positive for sufficiently large~$j$ (namely for
$j>-k_1$), which proves the desired convergence.

Finally, $\prec_j$ is different from $\prec_{_D}$ for all positive integers~$j$, 
since its smallest positive element is the conjugate of~$\sigma_2$ by 
$\sigma_2^{-j}\sigma_1$, and this is different from $\sigma_2$.
\label{cle}
\end{example}

\begin{figure}[htb]
\centerline{\input{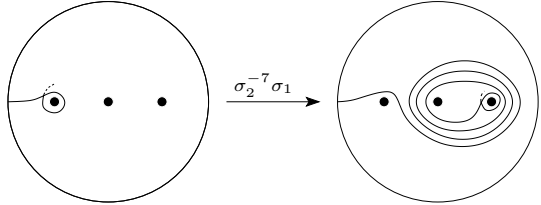}}
\vspace{-2mm}
\caption{The geodesic $\sigma_2^{-j}\sigma_1.\gamma$ deviates 
from~$\gamma$ to the left, but ``by as little as possible'':
it goes ``over'' the leftmost puncture (thus deviating to the left), 
but then it goes as much to the right as possible (it turns $j$ times 
counterclockwise around the two rightmost punctures, with~$j$ 
arbitrarily large)}
\label{F:DehorApprox}
\end{figure}

\begin{rmk}\rm The choice of the sequence $\sigma_2^{-j} \sigma_1$ 
is quite natural. Indeed, the sequence $\sigma_2^{-j}\sigma_1$ approaches 
the Conradian soul from above: for any element $c$ of the Conradian soul 
$\langle\sigma_2\rangle$ and for any $\sigma_1$-positive word~$w$ we have 
$c \prec_{_D} \sigma_2^{-j}\sigma_1 \prec_{_D} w$, where the second 
inequality is only true for sufficiently large~$j$.
In terms of geodesics, 
if~$\gamma$ is the standard geodesic inducing the Dehornoy-order
on~$B_3$ (c.f.\ Figure~\ref{F:D4geod}(a)), then the sequence
of geodesics $\sigma_2^{-j}\sigma_1.\gamma$ is eventually further to
the right than any prescribed geodesic $\beta.\gamma$, for $\beta$ a
$\sigma_1$-positive braid -- see Figure~\ref{F:DehorApprox}.
\end{rmk}

Example \ref{cle} can be generalized as follows.

\medskip

\begin{thm}\label{T:ConradZApproxConj}
Every Nielsen-Thurston ordering $\prec$ of finite type whose Conradian 
soul is a copy of $B_2$ can be approximated by its conjugates.
\end{thm}

{\bf Proof.} If $\sigma$ is a positive generator of the Conradian soul, then $\sigma$ 
corresponds to a left half-twist. According to \cite[Proposition 4.9]{navas-dynamics} 
and its proof (see also \cite[Theorem 2.9]{CNR} and its proof), there exists a 
sequence $\prec_j = g_j(\prec)$ of conjugates of $\prec$  which converges either 
to~$\prec$ or to~$\prec'$, where~$\prec'$ is the ordering whose set of positive 
elements is the union of $\{\sigma^k: k < 0\}$ and 
$\{g: g \succ 1, g \notin \langle \sigma \rangle \}$. If the limit point was $\prec'$, 
then the conjugates $g_j^{-1} \sigma g_j$ would be negative for $j$ large enough.  
This contradicts the fact that the Nielsen-Thurston ordering~$\prec$ has the 
subword property~S.

\begin{rmk}\label{R:PermGens}
\rm The argument of the proof above actually shows the stronger result that for any given 
Nielsen-Thurston ordering of finite type $\prec$, there exists another finite-type 
Nielsen-Thurston ordering $\prec_1$, obtained from the given one simply by permuting 
the lexicographic order of the Conradian soul, which is approximated 
by conjugates of itself (compare \cite[Proposition 4.9]{navas-dynamics} and 
\cite[Corollary 2.10]{CNR}). 
\end{rmk}

The previous results show that the only Nielsen-Thurston orderings which may fail to be 
approximated by their conjugates are those of finite type whose Conradian soul is 
higher-rank Abelian. 

\begin{figure}[ht]
\centerline{\epsfig{file=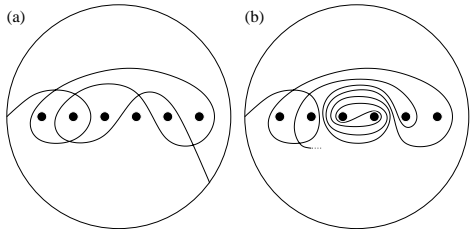,height=3cm}}
\caption{(a) A geodesic whose associated ordering of~$B_6$ 
seems not to be a limit point of Nielsen-Thurston orderings. 
Note the convex subgroups $\langle\sigma_5\rangle\subseteq 
\langle\sigma_3,\sigma_5\rangle\subseteq 
\langle\sigma_1,\sigma_3,\sigma_5\rangle$.
(b) Conjugating by $\sigma_3^{-N}\sigma_4$ (a sequence
which becomes eventually smaller than any positive
element outside the Conradian soul) 
yields an approximation of a different ordering with
convex jumps $\langle\sigma_1\rangle\subseteq 
\langle\sigma_1,\sigma_5\rangle\subseteq 
\langle\sigma_1,\sigma_3,\sigma_5\rangle$.}
\label{F:Counterex}
\end{figure}

\begin{example}\rm Consider the ordering~$\prec$ induced by the geodesic shown
in Figure~\ref{F:Counterex}(a). The Conradian soul of this ordering is 
$\langle\sigma_1,\sigma_3,\sigma_5\rangle=\Z^3$. 
By Theorem~\ref{T:NTnotIsol}, this order can be approximated by other
orders, but we do not know whether it can be approximated by its conjugates.
Imitating Example \ref{cle} and the proof of Theorem~\ref{T:ConradZApproxConj},
we could try to conjugate the order by elements which are as small as possible while
being larger than all elements of the Conradian soul. This, however, yields   
an approximation not of the desired ordering~$\prec$, but of one where the 
lexicographic order of $\mathbb{Z}^3$ has been permuted.
\end{example}

{\bf Question.} Can the ordering induced by the geodesic of 
Figure~\ref{F:Counterex}(a) be approximated by its conjugates? 
If not, can it be approximated by a sequence of Nielsen-Thurston 
orderings? The calculation above suggests that the answer to
both questions may be ``no''.

\bigskip


\noindent{\bf Acknowledgments.} The first-named author was funded 
by PBCT-CONICYT Research Network on Low Dimensional Dynamics.


\begin{small}

\end{small}

\end{document}